\documentclass[a4paper,  leqno,  myheadings,  11pt,  twoside]{article}

\usepackage{amsmath,  amsthm,  amssymb,  amscd,  amsxtra, graphicx}
\usepackage{latexsym,  amsfonts,enumerate}

\usepackage{color}
\definecolor{mycolor}{rgb}{0,1,0}
\usepackage[bookmarks, colorlinks, linkcolor=mycolor]{hyperref}

\def\T{Teich\-m\"ul\-ler }

\def\q1s{Q^1(S)}

 \makeatletter
\@addtoreset{equation}{section}

\begin{document}

\title{\bf{A binary
infinitesimal form of  \T metric
 }}
\author{GUOWU YAO
}
 \date{May 14,  2006}
\maketitle
\begin{abstract}\noindent
Let $S$ be a  Riemann surface of analytic finite type or the unit
disk in the complex plane. Let $[\mu]$ denote  the Teich\-m\"ul\-ler
equivalence classes of Beltrami differentials $\mu $. We apply the
Fundamental Inequalities to obtain a binary infinitesimal form of
Teich\-m\"ul\-ler metric. Using this form,  we define
``\emph{angle}" between two geodesics originating from a point  and
conjecture that the sum of the angles of a triangle in $T(S)$ should
be less than $\pi$ if $S$ is of analytic finite type. As a
consequence, the well-known necessary condition for two geodesics
coinciding is derived immediately.
\end{abstract}
\renewcommand{\thefootnote}{}

 \footnote{The   author was
supported by the National Natural Science Foundation of China (Grant
No. 10401036) and a Foundation for the Author of National Excellent
Doctoral Dissertation (Grant No. 200518) of PR China.}

\renewcommand\refname{\centerline{\Large{R}\normalsize{EFERENCES}}}

 \noindent
Department of Mathematical Sciences\\
  Tsinghua University\\Beijing,  100084,  People's Republic of
  China \\
  e-mail: \texttt{gwyao@math.tsinghua.edu.cn}
  \end{document}